\documentclass[12pt]{amsart}
\usepackage{amsfonts, amssymb}

\def\Cliff{{\rm Cliff}}

\newtheorem{thm}{Theorem}[section]
\newtheorem{prop}[thm]{Proposition}
\newtheorem{lemma}[thm]{Lemma}
\newtheorem{cor}[thm]{Corollary}
\def\proof{{\sc Proof. }}

\def\R{\mathbb{R}}

\def\CP{\mathbb{C}{\rm P}}
\def\C{\mathbb{C}}
\def\H{\mathbb{H}}
\def\N{\mathbb{N}}
\def\Box{\square}
\def\<{\langle}
\def\>{\rangle}

\def\Im{{\rm Im}}
\def\le{\leqslant}

\def\d{\partial}
\def\Nc{{\mathcal N}}

\def\rk{{\rm rk}}

\title{Circles and Quadratic Maps Between Spheres}
\author[V. Timorin]{Vladlen Timorin}
\thanks{Partially supported by CRDF RM1-2086}
\address{Institute for Mathematical Sciences,
Stony Brook University, Stony Brook NY 11794-3660, USA}
\email{<timorin@math.sunysb.edu>}
\keywords{line, circle, quadratic map between spheres, normed pairing}
\subjclass[2000]{15A63}

\begin{document}

\begin{abstract}
Consider an analytic map of a neighborhood of 0 in a vector space to
a Euclidean space.
Suppose that this map takes all germs of lines passing through 0 to germs of circles.
Such a map is called rounding.
We introduce a natural equivalence relation on roundings and
prove that any rounding, whose differential at 0 has rank at least 2, is
equivalent to a fractional quadratic rounding.
A fractional quadratic map is just the ratio of a quadratic map and a quadratic polynomial.
We also show that any rounding gives rise to a quadratic map between spheres.
The known results on quadratic maps between spheres have
some interesting implications concerning roundings.
\end{abstract}

\maketitle

\section*{Introduction}

By a {\em circle} in a Euclidean space we mean a round Euclidean circle, or a
straight line, or a point.
By a {\em vector line} in a vector space we mean a line passing through the origin.

Fix a vector space structure on $\R^m$ and a Euclidean structure on $\R^n$.
Consider a germ of an analytic map $\Phi:(\R^m,0)\to (\R^n,0)$.
We say that $\Phi$ is a {\em rounding} if it takes all germs of vector lines to germs of circles.
By the {\em rank of a rounding} we mean the rank of its first differential at 0.

\bigskip
\begin{tabular}{|c|}
\hline
Throughout this paper, we will always assume that\\ the rank of
a rounding is at least 2.\\
\hline
\end{tabular}
\bigskip

Roundings of ranks 1 and 0 are also interesting but
their study requires different methods.

We say that two roundings $\Phi_1$ and $\Phi_2$ are {\em equivalent} if
for any germ of vector line $l$ the germs $\Phi_1(l)$ and $\Phi_2(l)$
belong to the same circle.
We are mostly interested in description of roundings up to the equivalence.

Roundings mapping $(\R^2,0)$ to $(\R^2,0)$ were described by A. Khovanskii in \cite{Kh}.
The problem was motivated by nomography (see \cite{GKh}).
{\em Nomograms} are used for a graphical representation of functions.
They usually involve one or several families of curves that provide the means of reading data.
The most common types of nomograms are nomograms with aligned points and circular nomograms.
How to transform one type of nomograms into the other?
This problem was posed by G.S. Khovanskii \cite{GKh} in 1970-s.
A description of all roundings in dimension 2 was a major part in the solution of this problem.

All roundings mapping $(\R^2,0)$ to $(\R^2,0)$ are equivalent to M\"obius transformations.
In other words, there is a point $p\ne 0$ in $\R^2$ with the following property.
For any germ of vector line $l$, the image $\Phi(l)$ belongs to a circle passing through $p$.
Izadi \cite{Iz} proved the same result for roundings mapping $(\R^3,0)$ to
$(\R^3,0)$ and having an invertible first differential at 0.

It turns out that in dimension 4, an analogous statement is wrong.
The simplest counterexample is as follows.
Let us identify $\R^4$ with $\C^2$ so that multiplication by $i$ is an orthogonal operator,
and embed $\C^2$ to the complex projective plane $\CP^2$.
Consider the germ at 0 of a complex projective transformation $P:\CP^2\to\CP^2$
such that $P(0)=0$.
It takes germs of real vector lines to circles.
But it is not equivalent (as a rounding) to a M\"obius transformation.

Nevertheless, there is a simple description of all roundings mapping $(\R^4,0)$ to $(\R^4,0)$
(see \cite{Tim}).
There are 2 natural roundings from
$(\Im(\H)\times\H,0)$ to $(\H,0)$ where $\H$ is the skew-field of
quaternions and $\Im(\H)$ is the set of all purely imaginary
quaternions. The first rounding sends $(x,y)$ to $(1+x)^{-1}y$ and
the second to $y(1+x)^{-1}$. Any rounding mapping $(\R^4,0)$ to
$(\R^4,0)$ is equivalent to a rounding obtained from one of
these 2 roundings by composing it with an $\R$-linear map
$\R^4\to\Im(\H)\times\H$.

The main theorem of this article is the following:

\begin{thm}
\label{rounding}
Any rounding is equivalent to some fractional quadratic
rounding, i.e. a rounding of the form $P/Q$ where $P$ is a map given in
coordinates by quadratic polynomials, and $Q$ is a quadratic polynomial.
\end{thm}

In Section 3, we give a definition of a degenerate rounding and show that
any such rounding factors through a linear projection to a smaller space.
A nondegenerate quadratic rounding gives rise to a quadratic map between
spheres, i.e. a map between Euclidean spaces that is given in coordinates
by quadratic polynomials and that takes the unit sphere to the unit sphere.
Thus the description of nondegenerate roundings is reduced to
description of quadratic maps between spheres.

Results of \cite{Yiu1,Yiu2} lead to some interesting consequences regarding
roundings, including the following:

\begin{thm}
\label{kappa}
There exists a nondegenerate rounding $\Phi:(\R^m,0)\to (\R^n,0)$
if and only if $n\ge \kappa(m)$, where $\kappa$ is an explicit function
introduced in \cite{Yiu2}.
\end{thm}

In Section 5, we will recall the definition of $\kappa$.

The article is organized as follows.
We introduce a complexification of the notion of circle in Section 1.
With the help of it, in Section 2, we establish a crucial algebraic property
of the Taylor expansion of a rounding.
Section 3 contains the proof of Theorem \ref{rounding}.
Closely related with fractional quadratic roundings are normed pairings
and quadratic maps between spheres which are briefly discussed in Sections 4 and 5.
Results of Yiu \cite{Yiu2} on quadratic maps between spheres are used in Section 5 to
describe all possible dimensions $m$ and $n$ for which there is a
nondegenerate rounding mapping $(\R^m,0)$ to $(\R^n,0)$.

I am grateful to A. Khovanskii for useful discussions and to the anonymous
referee for useful comments and suggestions.

\section{Complex circles}

Let $\R^n$ be a Euclidean space with the Euclidean inner product $\<\cdot,\cdot\>$.
Consider the complexification $\C^n$ of $\R^n$ and extend the inner product to it by bilinearity.
The extended complex bilinear inner product will be also denoted by $\<\cdot,\cdot\>$.
Note that this inner product is not Hermitian!
We will never use Hermitian inner products.

A {\em complex circle} in $\C^n$ is either a line, or the intersection of an affine
2-plane (2-dimensional plane) with a {\em complex sphere}, i.e. a quadratic hypersurface
$\{x\in\C^n|\ \<x,x\>=\<a,x\>+b\}$ where $a\in\C^n$ and $b\in\C$.
Note that the complex sphere is a complex algebraic hypersurface, not a round sphere!

\begin{prop}
\label{Ccirc1}
A complex circle is a smooth curve, or a pair of intersecting lines, or a 2-plane.
\end{prop}

\proof
Consider a complex circle $C=P\cap S$, where $P$ is a 2-plane and $S$ is a complex sphere.
If $P$ does not belong to tangent hyperplanes to $S$ at points of $C$, then the
intersection of $P$ and $S$ is transverse, hence $C$ is smooth.
Now suppose that $P$ belongs to the tangent hyperplane to $S$ at a point $x\in C$.
The intersection $S\cap T_xS$ is a quadratic cone centered at the point $x$, which lies in $P$.
If we now intersect this quadratic cone with
$P$, then we obtain either the whole plane $P$ or a pair of intersecting
(possibly coincident) lines from $P$.
$\Box$

The {\em null-cone} (isotropic cone) $\{x\in\C^n|\ \<x,x\>=0\}$ will be denoted by $\Nc$.

\begin{prop}
\label{curve}
Consider a germ of holomorphic curve $\gamma:(\C,0)\to (\C^n,0)$, whose
image lies in a complex circle.
Suppose that $\gamma'(0)$ is nonzero and belongs to $\Nc$.
Then the linear span of the image of $\gamma$ belongs entirely to $\Nc$.
\end{prop}

\proof
Consider a complex circle $C$ containing the image of $\gamma$.
If $C$ is a line or a pair of lines, then the statement is obvious
(the curve $\gamma$ cannot switch from one line to the other).

Suppose that $C$ is the intersection of some plane and a complex sphere $S$.
Since $S$ contains the origin, its equation has the form
$$
\<x,x\>=\<a,x\>,\quad a\in\C^n.
$$
The curve $\gamma$ belongs to $S$.
Therefore, $\<\gamma,\gamma\>=\<a,\gamma\>$ identically.
Differentiating this relation at 0, we obtain that $\<a,v\>=0$ where $v=\gamma'(0)$.
It follows that the whole vector line spanned by $v$ belongs to $S$ and hence
to the circle $C$ containing the image of $\gamma$.

We see that $C$ is a 2-plane $P$.
This 2-plane $P$ must belong entirely to $\Nc$.
Indeed, if the restriction of the inner product to $P$ were nontrivial,
then in some coordinates $(x_1,x_2)$ on $P$,
it would be given by $x_1^2+x_2^2$ or $x_1^2$.
In both cases, a point with sufficiently large real $x_1$ and a real $x_2$
would not be in $S$.
$\Box$

\begin{prop}
\label{Ccirc2}
A set $X\subset\R^n$ lies in a complex circle if and only if it lies
in a real circle.
\end{prop}

\proof
Suppose that $X\subset\R^n$ lies in a complex circle $C$.
If $C$ is a complex line, then $C\cap\R^n$ is a real line containing $X$.
Otherwise, $C=P\cap S$, where $P$ is a complex 2-plane and $S$ is a complex sphere.
Then $X$ belongs to the real 2-plane $P\cap\R^n$.
Let $S$ be given by the equation
$$
\<x,x\>=\<a,x\>+b,\quad a\in\C^n,\ b\in\C.
$$
The real part of this equation defines a real sphere in $\R^n$.
Thus $X$ lies in the intersection of a real 2-plane with a real sphere.
But this intersection is a real circle.

In the other direction, the statement is obvious.
$\Box$

\begin{prop}
\label{Ccirc3}
A set $X\subset\C^n$ lies in a complex circle passing through 0 if and
only if the set of points $(x,\<x,x\>)$, $x\in X$, spans at most 2-dimensional
vector subspace of $\C^{n+1}$.
\end{prop}

\proof
Assume that $X$ lies in a complex circle passing through 0.
If it lies in a vector line $L$, then all vectors of the form $(x,\<x,x\>)$ are linear
combinations of $(a,0)$ and $(0,1)$, where $a$ is any point of $L$.

If $X$ does not belong to a single line, then it lies in some sphere $S$ containing the origin.
The restriction of the inner product $\<\cdot,\cdot\>$ to $S$ equals to some linear function
restricted to $S$.
In an orthonormal coordinate system, we have
$$
\<x,x\>=\lambda_1x_1+\dots+\lambda_nx_n
$$
for any point $x$ from $X$, where $x_i$ are the coordinates of $x$.
The coefficients $\lambda_i$ are independent of $x$.
Therefore, the system of vectors $(x,\<x,x\>)_{x\in X}$ has the same rank as
the set $X$.
But a circle is a plane curve, therefore, $X$ has rank 2.

The proof in the opposite direction is a simple reversion of the argument given above.
$\Box$

\begin{prop}
\label{circder}
Let $\Phi$ be an analytic map of a neighborhood of 0 in $\C^m$ to
$\C^n$. Take a vector $v\in\R^m$. The $\Phi$-image of the germ of
the real line spanned by $v$ lies in a complex circle if and only if the set
of vectors $(\d_v^k\Phi(0),\d_v^k\<\Phi,\Phi\>(0))_{k\in\N}$
spans at most 2-dimensional vector subspace of $\C^{n+1}$. Here $\d_v$
denotes the differentiation along the constant vector field $v$.
\end{prop}

\proof
By the Taylor formula, the vector subspace spanned by the image
of an analytic map coincides with the vector subspace spanned by all its
derivatives.
Consider the analytic map of $\R$ to $\C^n\times\C$ taking $t\in\R$ to
$(\Phi(tv),\<\Phi(tv),\Phi(tv)\>)$.
The vector space spanned by the image of this map is exactly
the linear span of the vectors $(\d_v^k\Phi(0),\d_v^k\<\Phi,\Phi\>(0))_{k\in\N}$.

It remains to use Proposition \ref{Ccirc3}.
$\Box$

\section{Roundings}

In this Section, we will give an algebraic constrain on the 2-jet of a rounding.

\begin{lemma}
Suppose that the linear span of vectors $y_1,y_2\in\C^n$ belongs to the
null-cone $\Nc$. Then we have
$$
\<y_1,y_1\>=\<y_2,y_2\>=\<y_1,y_2\>=0.
$$
\end{lemma}

\proof
Indeed, for any complex number $\lambda\in\C$, we have $y_1+\lambda y_2\in\Nc$.
Equating coefficients with all powers of $\lambda$ in the equation
$\<y_1+\lambda y_2,y_1+\lambda y_2\>=0$ to 0, we obtain the desired result.
$\Box$

Let $\Phi:(\R^m,0)\to(\R^n,0)$ be a rounding.
Denote by $A$ its first differential at 0.
In other words, $A:\R^m\to\R^n$ is a linear map such that $\Phi(x)=A(x)+\alpha(x)$,
where $|\alpha(x)|=o(|x|)$ as $x\to 0$.
Denote by $\<A,A\>$ the quadratic form on $\R^m$ whose value on a vector $x\in\R^m$
is $\<A(x),A(x)\>$.
Analogously, we define analytic functions $\<\Phi,\Phi\>$ and $\<\Phi,A\>$ on $\R^m$:
$$
\<\Phi,\Phi\>:x\mapsto\<\Phi(x),\Phi(x)\>,\quad \<\Phi,A\>:x\mapsto\<\Phi(x),A(x)\>.
$$

\begin{thm}
\label{div1}
The analytic functions
$\<A,\Phi\>$ and $\<\Phi,\Phi\>$ are divisible by $\<A,A\>$ in
the class of formal power series.
\end{thm}

This theorem gives conditions on an analytic map $\Phi$
that are necessary for $\Phi$ being a rounding.
We will only need restrictions that are imposed on the 2-jet of $\Phi$
by Theorem \ref{div1}.
It turns out (see Section 3) that these restrictions are sufficient to
guarantee that there exists a rounding with a given 2-jet.

\proof
Since $\Phi$ is analytic, it admits an analytic continuation to a neighborhood
of 0 in $\C^m$.

By Proposition \ref{circder}, we have that for any $v\in\R^m$ the set
of vectors $(\d_v^k\Phi(0),\d_v^k\<\Phi,\Phi\>(0))_{k\in\N}$
spans at most 2-dimensional subspace of $\R^{n+1}$. But this is
an algebraic condition on the coefficients of $\Phi$. Thus it holds
for all $v\in\C^m$. Using Proposition \ref{circder} again, we conclude
that the image of the germ at 0 of any complex vector line from $\C^m$
lies in some complex circle.

Suppose that a point $x\in\C^m$ is such that $\<A(x),A(x)\>=0$ but $A(x)\ne 0$.
Consider the holomorphic curve $\gamma:\C\to\C^n$ defined by the formula
$\gamma(\tau)=\Phi(\tau x)$ for all $\tau\in\C$.
The image of this holomorphic curve is the image of the complex line
in $\C^m$ spanned by $x$.
Therefore, the image of $\gamma$ lies in some complex circle.

Moreover, we have $\gamma'(0)=A(x)\ne 0$ and $\<\gamma'(0),\gamma'(0)\>
=\<A(x),A(x)\>=0$.
The latter means that $\gamma'(0)\in\Nc$, where $\Nc$ is the null-cone.
By Proposition \ref{curve}, the linear span $M$ of the image of
$\gamma$ belongs to $\Nc$.

Denote by $\Phi_k$ the power series of $\Phi$ at 0 truncated at degree $k$.
In other words, $\Phi(x)=\Phi_k(x)+\alpha_k(x)$, where
$|\alpha_k(x)|=o(|x|^k)$.
The vector subspace $M\subseteq\C^n$ contains the image of $\gamma$.
Therefore, it contains all derivatives of $\gamma$ at 0.
Then it must also contain
$$
\Phi_k(x)=\gamma'(0)+\frac{\gamma''(0)}2+\cdots+\frac{\gamma^{(k)(0)}}{k!}.
$$

Since $\Phi_k(x)$ and $A(x)=\Phi_1(x)$ belong to $M$, the linear
span of $\Phi_k(x)$ and $A(x)$ belongs to $\Nc$.
By the preceding lemma, we have
$$
\<A(x),\Phi_k(x)\>=\<\Phi_k(x),\Phi_k(x)\>=0.
$$
This equality holds for all $x$ satisfying the condition
$\<A(x),A(x)\>=0$.

By definition of a rounding, the rank of $A$ is et least 2.
Now we have two cases: $\rk(A)>2$ and $\rk(A)=2$.

{\em Case 1.}
In the first case, the equation $\<A,A\>=0$ defines an
irreducible hypersurface $\Gamma$.
Any polynomial vanishing on $\Gamma$ is divisible by $\<A,A\>$, the equation of $\Gamma$.
Therefore, for any $k$, the polynomials $\<A,\Phi_k\>$ and $\<\Phi_k,\Phi_k\>$ are
divisible by $\<A,A\>$.
The theorem now follows.

{\em Case 2.}
In the second case, the equation $\<A,A\>=0$
splits into the product of two distinct linear factors: $\<A,A\>=l_1l_2$.
Here $l_1$ and $l_2$ are linear functionals with complex coefficients.
The equation $\<A,A\>=0$ defines the union of two hyperplanes: one is
given by the equation $l_1=0$ and the other is given by the equation $l_2=0$.
Both polynomials $\<A,\Phi_k\>$ and $\<\Phi_k,\Phi_k\>$ vanish on
the hyperplane $l_1=0$.
Therefore, both polynomials are divisible by $l_1$.
Analogously, both polynomials are divisible by $l_2$.
Thus they are divisible by $\<A,A\>=l_1l_2$.
$\Box$

Theorem \ref{div1} has two important corollaries.
Recall that a map $B:\R^m\to\R^n$ is a {\em quadratic map} if
for any constant vector field $v$ on $\R^m$ the third derivative
$\d_v^3B$ vanishes everywhere.

\begin{cor}
\label{div2}
Let $\Phi$ be a rounding.
Denote by $A$ its linear part and by $B$ its quadratic part.
In other words, $A$ is a linear map and $B$ is a homogeneous quadratic
map of $\R^m$ to $\R^n$ such that $\Phi(x)=A(x)+B(x)+C(x)$, where
$|C(x)|=o(|x|^2)$.
Then both polynomials $\<A,B\>$ and $\<B,B\>$ are divisible by $\<A,A\>$.
\end{cor}

\proof
We will use the notation of the proof of Theorem \ref{div1}.
We have $\Phi_2=A+B$.
Polynomials $\<A,\Phi_2\>$ and $\<\Phi_2,\Phi_2\>$ are divisible by $\<A,A\>$.
This follows from Theorem \ref{div1}, and even more directly, it follows from
the proof of it.
Since $\<A,A+B\>$ is divisible by $\<A,A\>$, the polynomial $\<A,B\>$ must also
be divisible by $\<A,A\>$.
Since $\<A+B,A+B\>=\<A,A\>+2\<A,B\>+\<B,B\>$ is divisible by $\<A,A\>$,
the polynomial $\<B,B\>$ must also be divisible by $\<A,A\>$.
$\Box$

\begin{cor}
\label{kerA}
Let $\Phi$ be a rounding with the first differential $A$. Then the kernel
of $A$ maps to 0 under $\Phi$.
\end{cor}

Now we are going to establish a criterion of the equivalence of two roundings.

\begin{lemma}
\label{2jet1}
If two roundings have the same 2-jets, then they are equivalent.
\end{lemma}

\proof
Let $\Phi,\Psi:(\R^m,0)\to (\R^n,0)$ be analytic roundings with
the same 2-jet $A+B$, where $A$ is the linear part and $B$ is the quadratic part.
It is clear that if a line from $\R^m$ does not belong to the kernel of
$A$, then it goes to the same circle under both maps.
Indeed, a circle is determined by the velocity (the first derivative) and
the acceleration (the second derivative) at 0 of any parameterized curve
lying in this circle and such that the velocity at 0 does not vanish.
Corollary \ref{kerA} concludes the proof.
$\Box$

We need the following technical fact:

\begin{lemma}
\label{paral}
Consider a linear map $A:\R^m\to\R^n$ of rank at least 2 and a quadratic
(or linear) homogeneous map $B:\R^m\to\R^n$ such that $B$ is everywhere parallel to $A$.
In other words, for any $x\in\R^m$, the vectors $A(x)$ and $B(x)$ are linearly dependent.
Then $B(x)=l(x)\cdot A(x)$ for some linear (or constant) function $l:\R^m\to\R$
and all $x\in\R^m$.
\end{lemma}

\proof
We will work out only the case when $B$ is a homogeneous quadratic map.
The linear case is only easier.
The function $l=B/A$ (mapping any vector $x\in\R^m$ to the real number $l(x)$ such that
$B(x)=l(x)A(x)$) is defined on the complement to the kernel of $A$.
The parallelogram equality for $B$ reads
$$
l(x+y)A(x+y)+l(x-y)A(x-y)=2(l(x)A(x)+l(y)A(y))
$$
where $x$ and $y$ are vectors from $\R^m$ such that none of the
vectors $x$, $y$, $x+y$ and $x-y$ lies in the kernel of $A$.

Suppose that $A(x)$ and $A(y)$ are linearly
independent.
Equating the coefficients with $A(x)$ in the parallelogram
equality, we obtain $l(x+y)+l(x-y)=2l(x)$.
Set $u=x+y$ and $v=x-y$.
Then $l(u+v)=l(u)+l(v)$.
This holds for almost all pairs $(u,v)$ of vectors from $\R^m$.
Therefore, $l$ extends to a linear function.
$\Box$

\begin{prop}
\label{2jet2}
Consider two roundings $\Phi=A+B+\dots$ and $\Phi'=A'+B'+\dots$,
where $A$ and $A'$ are linear maps, $B$ and $B'$ are homogeneous
quadratic maps, and dots denote terms of order 3 and higher.
The roundings $\Phi$ and $\Phi'$ are equivalent if and
only if $A'(x)=\lambda A(x)$ and $B'(x)=\lambda^2 B(x)+l(x)A(x)$, where
$\lambda$ is a real number, and $l:\R^m\to\R$ is a linear functional.
\end{prop}

\proof
First suppose that $A'$ and $B'$ are related to $A$ and $B$ as above.
Composing $\Phi$ with the local diffeomorphism $x\mapsto \lambda x+l(x)x$,
which preserves all germs of vector lines, we obtain a rounding equivalent
to $\Phi$ and having the 2-jet $A'+B'$.
By Lemma \ref{2jet1}, this rounding is equivalent to $\Phi'$.
Thus $\Phi$ and $\Phi'$ are equivalent.

Now suppose that for any germ of vector line in $\R^m$
the roundings $\Phi$ and $\Phi'$ map this germ to the same circle.
In particular, the vectors $A(x)$ and $A'(x)$ are linearly dependent for
all $x\in\R^m$.
Indeed, both vectors $A(x)$ and $A'(x)$ are tangent at 0 to the circle
containing the images of the vector line spanned by $x$ under the maps
$\Phi$ and $\Phi'$.
By Lemma \ref{paral}, we have $A'(x)=\lambda A(x)$ for a real number $\lambda$
independent of $x$.

The orthogonal projection of $B'(x)$ to the orthogonal complement of
$A(x)$ is the same as that of $\lambda^2 B(x)$.
This follows from the fact that the images of the line spanned by $x$
under the maps $\Phi$ and $\Phi'$ lie in the same circle, in particular,
they have the same curvature at 0.
Thus the vectors $B'(x)-\lambda^2 B(x)$ and $A(x)$ are linearly dependent
for all $x\in\R^m$.
The Proposition now follows from Lemma \ref{paral}.
$\Box$

\section{Fractional quadratic maps}

Let $F$ be a quadratic map of $\R^m$ to $\R^n$, and $Q$ a quadratic function
on $\R^m$ ($F$ and $Q$ are not necessarily homogeneous).
The map $F/Q:x\mapsto F(x)/Q(x)$ is called {\em a fractional quadratic map}.
It is defined on the complement to the zero level of $Q$.
Nonetheless, it will be referred to as a fractional quadratic map of $\R^m$ to $\R^n$.

A map from an open subset $U$ of $\R^m$ to $\R^n$ is said to take all lines
to circles if any germ of line contained in $U$ gets mapped to a germ of circle
under this map.

\begin{prop}
\label{fq}
Let $F:\R^m\to\R^n$ be a quadratic map such that the polynomial
$\<F,F\>:x\mapsto\<F(x),F(x)\>$ is divisible by some quadratic function $Q$.
Then the fractional quadratic map $F/Q$ takes all lines to circles.
\end{prop}

\proof
Introduce an orthonormal basis in $\R^n$.
Denote the components of $F(x)$ with respect to this basis by $F_1(x),\dots,F_n(x)$.
Thus $F_1,\dots,F_n$ are quadratic polynomials on $\R^m$, not necessarily
homogeneous.

Take an arbitrary line $L\subset\R^m$ not lying in the zero level of $Q$.
The functions $F_1,F_2,\dots,F_n$ and $Q$ restricted to $L$ span a subspace
$V$ of the space of all quadratic polynomials on $L$.
We have a natural map of the space of all affine functions on $\R^n$ to the space $V$.
Namely, an affine function $f$ on $\R^n$ gets mapped to the function
$x\mapsto Q(x)f(F(x)/Q(x)))$, which is a quadratic function belonging to the space $V$.
Clearly, this map is onto.
The kernel of this map, which consists of all affine functions vanishing
identically on the image of $L$ under the map $F/Q$, has dimension $n+1-\dim(V)$.

If the subspace $V$ is one-dimensional, then $(F/Q)(L)$ is just a point.
If $V$ is 2-dimensional, then $(F/Q)(L)$ belongs to a line.
Finally, if $V$ is 3-dimensional
(i.e. it contains all quadratic polynomials on $L$)
then $(F/Q)(L)$ is a plane curve.

Consider the last case.
The ratio $\<F,F\>/Q$ is a quadratic polynomial.
Its restriction to $L$ belongs to $V$.
Therefore, $\<F,F\>/Q$ is a linear combination of functions $F_i$ and $Q$ on $L$,
and $\<F,F\>/Q^2$ is a linear combination of $F_i/Q$ and $1$ on $L$.
This means that under restriction to the image of $L$,
the square of the Euclidean norm equals to some affine function.
Hence the image of $L$ lies in some sphere.
A plane curve that lies in a sphere is necessarily a circle.
$\Box$

The main result of this section is as follows:

\begin{thm}
\label{r1}
Let $A:\R^m\to\R^n$ be a linear homogeneous map and $B:\R^m\to\R^n$ a
quadratic homogeneous map such that both polynomials
$\<A,B\>$ and $\<B,B\>$ are divisible by $\<A,A\>$.
Then the map
$$
\Psi=\frac{A+B-2pA}{1-2p+q},\qquad p=\frac{\<A,B\>}{\<A,A\>},\ q=\frac{\<B,B\>}{\<A,A\>},
$$
is a fractional quadratic rounding whose 2-jet at 0 is $A+B$. Moreover,
$\Psi$ takes all lines to circles, not only those passing through 0.
\end{thm}

\proof
By our assumptions, $p$ is a linear polynomial, and $q$ is a quadratic polynomial.

It is readily seen that $\Psi$ has 2-jet $A+B$.
We also claim that $\Psi$ rounds all germs of lines
(i.e., maps them to germs of circles).
Indeed, this follows from Proposition \ref{fq}, since the square of the norm of
$A+(B-2pA)$ equals to $(1-2p+q)\<A,A\>$.
$\Box$

\begin{cor}
\label{r2}
Suppose that a rounding $\Phi$ has 2-jet $A+B$.
In other words, $A:\R^m\to\R^n$ is a linear map, and $B:\R^m\to\R^n$
is a homogeneous quadratic map such that $\Phi(x)=A(x)+B(x)+\cdots$,
where dots denote terms of order 3 and higher.
Then $\Phi$ is equivalent to the fractional quadratic rounding
$$\frac{A+B-2pA}{1-2p+q}$$
where $p=\<A,B\>/\<A,A\>$ and $q=\<B,B\>/\<A,A\>$.
\end{cor}

This follows from Theorem \ref{r1} and Lemma \ref{2jet1}.

Consider a rounding $\Phi$ with 2-jet $A+B$, where $A$ is a linear map,
and $B$ is a homogeneous quadratic map.
The rounding $\Phi$ is said to be {\em degenerate} if
there exists a point $x_0\in\R^m$ such that $A(x_0)=0$ and $q(x_0)=p^2(x_0)$,
where $p$ and $q$ are as above.
This definition may look artificial but it is explained by the following

\begin{lemma}
\label{degen}
A rounding equivalent to a degenerate rounding is itself degenerate.
Any degenerate rounding $\Phi$ with 2-jet $A+B$ is equivalent to a rounding
$\Phi'=A'+B'+\cdots$ such that the kernels of $A'$ and $B'$ intersect
nontrivially.
\end{lemma}

\proof
Let $\Phi$ be a degenerate rounding with 2-jet $A+B$.
Assume that a rounding $\Phi'=A'+B'+\cdots$ is equivalent to $\Phi$.
Then by Lemma \ref{2jet2}, we have $A'=\lambda A$ and $B'=\lambda^2 B+lA$, where $\lambda$ is a
number and $l$ is a linear function. The polynomials $p$ and $q$ of
these 2 roundings are related as follows:
$$
p'=\lambda p+\frac l{\lambda},\quad q'=\lambda^2 q+2lp+
\frac{l^2}{\lambda^2}.
$$
If $A=0$ at some point, then $A'=0$ at the same point.
If $q=p^2$ at some point, then $q'=p'^2$.
This proves the first part of the Lemma.

To prove the second part, choose $\lambda=1$, $l=-p$.
Then $p'=0$, $q'=q-p^2$.
Now assume that at a point $x_0\in\R^m$, we have $A(x_0)=0$ and $q(x_0)=p^2(x_0)$.
Then $q'(x_0)=0$, hence $x_0$ is a zero of $B'$ of order greater than 1,
so it lies in the kernel of $B'$.
Thus the point $x_0$ lies in the kernels of both $A'$ and $B'$.
$\Box$

We can concentrate on nondegenerate roundings only, due to the following

\begin{prop}
Let $\Phi:(\R^m,0)\to (\R^n,0)$ be a degenerate rounding. Then there is
a projection $\pi$ from $\R^m$ to a smaller space $\R^k$, $k<m$, and a nondegenerate
fractional quadratic rounding
rounding $\Psi:(\R^k,0)\to (\R^n,0)$ such that $\Psi\circ\pi$ is equivalent
to $\Phi$.
\end{prop}

\proof
By Lemma \ref{degen} we can assume that the intersection $K$ of
the kernels of $A$ and $B$ is a nontrivial subspace of $\R^m$.
Denote by $\pi$ the natural projection of $\R^m$ to the quotient
$\R^k=\R^m/K$.
Then $A=\tilde A\circ\pi$ and $B=\tilde B\circ\pi$ where
the maps $\tilde A$ and $\tilde B$ from $\R^k$ to $\R^n$ are linear and
quadratic, respectively. Both polynomials $\<\tilde A,\tilde B\>$ and
$\<\tilde B,\tilde B\>$ are divisible by $\<\tilde A,\tilde A\>$.

By Theorem \ref{r1} there is a fractional quadratic rounding $\Psi:(\R^k,0)\to (\R^n,0)$
with the 2-jet $\tilde A+\tilde B$. Hence $\Psi\circ\pi$ is a
fractional quadratic rounding
with the 2-jet $A+B$. By Lemma \ref{2jet1} it is equivalent to $\Phi$.
$\Box$

\section{Normed pairings}

It turns out that roundings are closely related with normed pairings
and quadratic maps between spheres.
In this section and in the next one, we review these subjects.

Hurwitz in 1898 posed the following problem which is still unsolved:
find all relations of the form
$$
(x_1^2+\cdots+x_r^2)(y_1^2+\cdots+y_s^2)=z_1^2+\cdots+z_n^2
$$
where $z_1,\dots,z_n$ are bilinear functions of the variables
$x_1,\dots,x_r$ and $y_1,\dots,y_s$.
A relation of this form is represented by a bilinear map
$f:\R^r\times\R^s\to\R^n$ such that for all $x\in\R^r$ and $y\in\R^s$
we have $|f(x,y)|=|x|\cdot|y|$.
Such maps are called {\em normed pairings} of size $[r,s,n]$.

Many examples of normed pairings are known.
The most familiar are the multiplication laws for real numbers, complex numbers,
quaternions and octonions.
These are normed pairings of sizes $[1,1,1]$, $[2,2,2]$,
$[4,4,4]$ and $[8,8,8]$, respectively.
As Hurwitz proved in 1898 \cite{Hu1}, these are
the only possible normed pairings of size $[n,n,n]$ up to orthogonal
transformations.
Later on, he succeeded in describing all normed pairings of size $[r,n,n]$, see \cite{Hu2}.
The same result was independently obtained by Radon \cite{Ra}.

Recall that the Clifford algebra $\Cliff(k)$ is the associative algebra
generated over reals by $r$ elements $e_1,\dots,e_k$ satisfying the relations
$$
e_i^2=-1,\qquad e_ie_j+e_je_i=0\quad (i\ne j).
$$
A linear representation of $\Cliff(k)$ in a Euclidean space $\R^n$ is called
{\em compatible with the Euclidean structure} if all generators $e_i$ act
as orthogonal operators.
Any finite dimensional representation of a Clifford
algebra is compatible with a suitable Euclidean structure on the space of
representation.
The result of Hurwitz and Radon is as follows:

\begin{thm}
Suppose that $f:\R^r\times\R^n\to\R^n$ is a normed pairing.
Then there is a representation $\phi$ of $\Cliff(r-1)$ in
$\R^n$ compatible with the Euclidean structure and such that
$$
f(x,y)=\phi(x_0+x_1e_1+\cdots+x_{r-1}e_{r-1})A(y)
$$
where $A$ is a linear conformal transformation and $x_0,\dots,x_{r-1}$ are
coordinates of $x$ in some orthonormal basis.
\end{thm}

The largest $r$ for which there is a representation of $\Cliff(r-1)$ in
$\R^n$ is denoted by $\rho(n)$ and is called the {\em Hurwitz--Radon function}
of $n$.
The Hurwitz--Radon theorem implies that for any normed pairing of
size $[r,n,n]$, we have $r\le\rho(n)$.
The Hurwitz--Radon function can be computed explicitly due to the result of
\'E. Cartan \cite{Car} who classified all Clifford algebras and their representations in 1908
(the idea of this classification goes back to Clifford).
See also \cite{ABS}.
Suppose that $n=2^su$, where $u$ is odd.
If $s=4a+b$, $0\le b\le 3$, then $\rho(n)=8a+2^b$.

\begin{prop}
Let $f:\R^r\times\R^s\to\R^n$ be a normed pairing.
The map $f$ can be regarded as a quadratic map of $\R^r\oplus\R^s$ to $\R^n$.
Consider the quadratic form $Q$ on $\R^r\oplus\R^s$
defined by the formula $Q(x,y)=\<x,x\>$ for all $x\in\R^r$, $y\in\R^s$.
Then the degree 4 polynomial $\<f,f\>$ on $\R^r\oplus\R^s$ is divisible by $Q$,
hence $f/Q$ takes all lines to circles.
\end{prop}

This Proposition follows immediately from definitions and Proposition \ref{fq}.
We see that normed pairings provide many examples of maps taking all
lines to circles.

A bilinear map $f:\R^r\times\R^s\to\R^n$ is called {\em nonsingular}
if from $f(x,y)=0$ it follows that $x=0$ or $y=0$.
Clearly, any normed pairing is a nonsingular bilinear map.
Possible sizes of nonsingular bilinear maps are restricted by the following theorem
of Stiefel \cite{St} and Hopf \cite{H}.

\begin{thm}
\label{HS}
If there is a nonsingular bilinear map of size $[r,s,n]$, then the
binomial coefficient $\binom{n}{k}$ is even whenever $n-r<k<s$.
\end{thm}

This is a topological theorem.
It uses the ring structure in the cohomology of projective spaces.

\section{Quadratic maps between spheres}

By a {\em unit sphere} $S^n$ we mean the set of all vectors in
$\R^{n+1}$ of Euclidean length 1.
A map $f:S^m\to S^n$ between unit spheres is called {\em quadratic},
if it extends to a quadratic homogeneous map of $\R^{m+1}$ to $\R^{n+1}$.
This extension must satisfy the condition $\<f(x),f(x)\>=\<x,x\>^2$ for all $x\in\R^{m+1}$.
By a {\em great circle} in $S^m$ we mean a circle obtained as the intersection
of $S^m$ with a 2-dimensional vector subspace.
The following simple but very important statement is proved in \cite{Yiu1}:

\begin{prop}
Any quadratic map $f:S^m\to S^n$ takes great circles to circles.
\end{prop}

Consider a homogeneous quadratic map $F:\R^{m+1}\to\R^n$
such that $\<F,F\>=Q_1\cdot Q_2$, where $Q_1$ and $Q_2$ are quadratic forms.
Then the forms $Q_1$ and $Q_2$ both take nonnegative values only
or both take non-positive values only.
By changing the sign, if necessary, we can always arrange that
both forms take nonnegative values only.
Suppose that the kernels of $Q_1$ and $Q_2$ intersect trivially.
Then $Q_1+Q_2$ is a positive definite form on $\R^{m+1}$.
Equip $\R^{m+1}$ with the Euclidean structure such that the
square of the length with respect to this structure equals to the
quadratic form $Q_1+Q_2$.
In other words, $\<x,x\>=Q_1(x)+Q_2(x)$ for any $x\in\R^{m+1}$.

The next proposition is verified by a simple direct computation:

\begin{prop}
\label{s}
The quadratic map $(2F,Q_1-Q_2):\R^{m+1}\to\R^{n+1}$
takes the unit sphere $S^m\subset\R^{m+1}$ to the unit sphere $S^n\subset\R^{n+1}$.
\end{prop}

Most examples of quadratic maps between spheres come from normed pairings.
Consider a normed pairing $f:\R^r\times\R^s\to\R^n$.
Then the {\em Hopf map}
$$
H_f(x,y)=(2f(x,y),\<x,x\>-\<y,y\>)
$$
takes $S^{r+s-1}$ to $S^n$.

Yiu \cite{Yiu2} described all pairs of positive integers $m,n$ such that there
is a non-constant quadratic map of $S^m$ to $S^n$.
Namely, $n\ge\kappa(m)$, where the Yiu function $\kappa$ is defined recurrently as follows:
$$
\kappa(2^t+m)=\left\{\begin{array}{lr} 2^t,& 0\le m<\rho(2^t)\\
2^t+\kappa(m),& \rho(2^t)\le m<2^t\end{array}\right.
$$

Let $f:S^m\to S^n$ be any map between unit spheres.
Denote by $\phi$ a map of an open subset $U$ of $\R^m$ to $\R^n$ obtained as
the composition of
\begin{itemize}
\item an affine embedding of $U\subseteq\R^m$ into $\R^{m+1}-0$,
\item the radial projection of $\R^{m+1}-0$ to $S^m\subset\R^{m+1}$,
\item the map $f:S^m\to S^n$,
\item a stereographic projection of $S^n$ to some hyperplane $H$ in $\R^{n+1}$,
\item a Euclidean identification of $H$ with $\R^n$.
\end{itemize}
Then we say that the map $\phi$ {\em factors through} the map $f$.

\begin{thm}
\label{r3}
Any nondegenerate rounding $\Phi:(\R^m,0)\to (\R^n,0)$ admits an equivalent
rounding that factors through a quadratic map between spheres $S^m$ and $S^n$.
\end{thm}

\proof
By Corollary \ref{r2}, there is a rounding $\Psi$ equivalent to $\Phi$
that extends to a fractional quadratic map $F/Q$, where $F:\R^m\to\R^n$
is a (possibly inhomogeneous) quadratic map, $Q:\R^m\to\R$ is a
(possibly inhomogeneous) quadratic polynomial and
the degree 4 polynomial $\<F,F\>$ is divisible by $Q$.
Extend $\R^m$ to $\R^{m+1}$ by adding an extra coordinate $t$.
Let $\tilde F$ and $\tilde Q$ be the homogeneous quadratic map and
the homogeneous quadratic form, respectively, that restrict to $F$ and $Q$ on
the hyperplane $t=1$.

Consider the map $\tilde F/\tilde Q:\R^{m+1}\to \R^n$.
By Proposition \ref{fq}, it sends all lines to circles.
Now compose this map with the inverse stereographic projection of $\R^n$
to the unit sphere in $\R^{n+1}$.
We obtain the map
$$
f':\R^{m+1}\to\R^{n+1},\quad f'=\left(\frac{2\tilde F}{\tilde Q_1+\tilde Q_2},
\frac{\tilde Q_1-\tilde Q_2}{\tilde Q_1+\tilde Q_2}\right),
\quad f'(\R^{m+1})\subseteq S^n,
$$
where $\tilde Q_1$ and $\tilde Q_2$ are quadratic forms such that
$\tilde F=\tilde Q_1\tilde Q_2$.
Without loss of generality, we can assume that forms $Q_1$ and $Q_2$
are nonnegative.

If the quadratic form $\tilde Q_1+\tilde Q_2$ is nondegenerate,
then it is positive definite, i.e. it defines a Euclidean inner
product on $\R^{m+1}$.
The quadratic map $f=(2\tilde F,\tilde Q_1-\tilde Q_2)$
from Proposition \ref{s} takes the unit sphere $S^m$ in $\R^{m+1}$
to the unit sphere $S^n$ in $\R^{n+1}$ and coincides with $f'$ on $S^m$.
Thus $\Psi$ factors through $f$.

It remains to verify that $\tilde Q_1+\tilde Q_2$ is indeed nondegenerate.
Using the explicit construction of $F$ and $Q$, we can write
$$
\tilde Q_1+\tilde Q_2=t^2-2pt+q+\<A,A\>=
(q-p^2)+(p-t)^2+\<A,A\>.
$$
Here $A+B$ is the 2-jet of $\Phi$ (and of $\Psi$), $p=\<A,B\>/\<A,A\>$
and $q=\<B,B\>/\<A,A\>$.
If at some point $x_0\in\R^m$, the form $\tilde Q_1+\tilde Q_2$ vanishes,
then $A(x_0)=0$, $q(x_0)=p^2(x_0)$, and hence $\Psi$ is degenerate.
Contradiction.
$\Box$

In \cite{Tim2}, we found a simple condition on a rounding, which guarantees
that it factors through the Hopf map associated with a representation
of a Clifford algebra.

Theorem \ref{r3} combined with results of Yiu leads to Theorem \ref{kappa}.


\begin{thebibliography}{9}
\bibitem{Kh}
A.G. Khovanskii: {\em Rectification of circles}, Sib. Mat. Zh.,
{\bf 21} (1980), 221--226

\bibitem{GKh}
G.S. Khovanskii: {\em Foundations of Nomography}, ``Nauka'', Moscow, 1976
(Russian)

\bibitem{Iz}
F.A. Izadi: {\em Rectification of circles, spheres, and
classical geometries}, PhD thesis, University of Toronto, (2001)

\bibitem{Tim}
V.A. Timorin: {\em ``Rectification of circles and quaternions''},
Michigan Mathematical Journal, \textbf{51} (2003), 153--167

\bibitem{Yiu1}
P. Yiu: {\em Quadratic forms between spheres and the non-existence
of sums of squares formulae}, Math. Proc. Cambridge Philos. Soc.,
{\bf 100}, 493-504

\bibitem{Yiu2}
P. Yiu: {\em Quadratic forms between euclidean spheres},
Manuscripta Math., {\bf 83}, 171--181

\bibitem{Hu1}
A. Hurwitz: {\em \"Uber die Komposition der quadratischen Formen von
beliebig vielen Variabeln}, Nachr. Ges. Wiss. G\"ottingen,
(Math.-Phys. Kl.) (1898) 309--316. Reprinted in  Math. Werke I, 565--571

\bibitem{Hu2}
A. Hurwitz: {\em \"Uber die Komposition der quadratischen Formen},
Math. Ann. 88 (1923) 1--25. Reprinted in Math. Werke II, 641--666

\bibitem{Ra}
J. Radon: {\em Lineare Scharen orthogonaler Matrizen}, Abh. Math. Sem.
Univ. Hamburg 1 (1922), 1--14

\bibitem{Car}
E. Cartan: {\em Nombres complexes}, pp. 329-448 in J. Molk (red.):
Encyclop\'edie des sciences math\'ematiques, Tome I, Vol. 1, Fasc. 4,
art. 15 (1908). Reprinted in E. Cartan: {\em \OE uvres compl\`etes},
Partie II. Gauthier-Villars, Paris, 1953, pp. 107--246

\bibitem{ABS}
M.F. Atiyah, R. Bott, A. Shapiro: {\em Clifford modules}, Topology {\bf 3},
suppl. 1 (1964), 3--38. Reprinted in R. Bott: {\em Lectures on $K(X)$}.
Benjamin, New York, 1969, pp. 143--178. Reprinted in M. Atiyah:
{\em Collected Works}, Vol. 2. Clarendon Press, Oxford, 1988, pp. 301--336

\bibitem{St} E. Stiefel: {\em \"Uber Richtungsfelder in den projektiven
R\"aumen und einen Satz aus der reelen Algebra}, Comment. Math. Helv.
{\bf 13} (1941), 201-218

\bibitem{H} H. Hopf: {\em Ein topologischer Beitrag zur reelen Algebra},
Comment. Math. Helv., {\bf 13} (1940/41) 219--239

\bibitem{Sh}
D. B. Shapiro: {\em Compositions of Quadratic Forms}, de Gruyter Expositions
in Math., {\bf 33} (2000)

\bibitem{Tim2}
V. Timorin: {\em Circles and Clifford Algebras},
Funktsional. Anal. i Prilozhen. \textbf{38} (2004), No. 1, 56--64

\bibitem{L}
K.Y. Lam: {\em Some new results in composition of quadratic forms}, Invent.
Math. {\bf 79} (1985), 467--474
\end{thebibliography}
\end{document}